\documentclass[a4paper,12pt]{article}
\usepackage{amssymb}

\newcommand{\cL}{{\mathcal{L}}}

\newcommand{\R}{\mathbb{R}}

\newcommand{\bS}{\mathbb{S}}

\newtheorem{theorem}{Theorem}[section]

\newtheorem{lemma}{Lemma}[section]

\newtheorem{remark}{Remark}[section]
\newcommand{\bremark}{\begin{remark} \em}
\newcommand{\eremark}{\end{remark} }

 \usepackage{cite}
 \usepackage{color}

\begin{document}

\begin{center}{\bf  \large   Qualitative properties for elliptic  problems\\[2mm]

 with CKN operators }\medskip
 \bigskip
 \medskip

{\small Huyuan Chen\footnote{chenhuyuan@yeah.net} \quad Yishan Zheng\footnote{zyszyszys15@126.com}  }
 \bigskip

 {\small
 Department of Mathematics, Jiangxi Normal University,\\
Nanchang, Jiangxi 330022, PR  China
\medskip

}

\bigskip

\begin{abstract}
The purpose of this paper is to study  basic property of
 the   operator
 $$\cL_{\mu_1,\mu_2} u=-\Delta   +\frac{\mu_1 }{|x|^2}x\cdot\nabla +\frac{\mu_2 }{|x|^2},$$ which generates at the origin due to the critical gradient and  the Hardy term, where $\mu_1,\mu_2$ are free parameters.  This operator arises from the critical Caffarelli-Kohn-Nirenberg inequality.   We analyze the fundamental solutions in a weighted distributional identity and   obtain the Liouville theorem for the Lane-Emden equation with that operator,
 by using the classification of isolated singular solutions of  the related Poisson problem in a bounded domain $\Omega \subset \R^N$ ($N \geq 2$) containing the origin.

\end{abstract}

\end{center}
   \vspace{1mm}


\vspace{2mm}

\setcounter{equation}{0}
\section{Introduction}
The well-known Caffarelli-Kohn-Nirenberg inequality (CKN inequality for short) proposed in \cite{CKN} states as following
$$
\Big(\int_{\R^N} |x|^{-b(p+1)} |u|^{p+1}dx\Big)^{\frac{2}{p+1}}\leq C_{a,b,N} \int_{\R^N}   |x|^{-2a} |\nabla u|^2 dx,
$$
where $N\geq 2$,
$$-\infty<a<\frac{N-2}{2},\quad a\leq b\leq a+1\quad{\rm and}\quad p=\frac{N+2(1+a-b)}{N-2(1+a-b)}.$$
A critical CKN inequality with $b=a+1$  (see \cite{CW,ACP,WW})   is
\begin{equation}\label{CNK-in10}
\int_{\R^N}|x|^{-2a}     |\nabla u|^2 dx \geq\Big(\frac{N-2-2a}{2}\Big)^2 \int_{\R^N} |x|^{-2(a+1)} |u|^2dx,
 \end{equation}
which, for $a=0 $  and $N\geq 3$,  reduces to the classical Hardy inequality
\begin{equation}\label{CNK-in1}
 \int_{\R^N}  |\nabla u|^2dx \geq  \frac{(N-2)^2}{4}\int_{\R^N}    \frac{|u|^2}{|x|^2}  dx.
 \end{equation}
 It is known that  the related  elliptic operator arising form (\ref{CNK-in10})
\begin{eqnarray*}
\cL_au&=&-{\rm div}(|x|^{-2a}\nabla u)-\frac{\big(N-2-2a\big)^2}{4}\frac{u}{|x|^{2(a+1)}}
\\&=&|x|^{-2a}\Big(-\Delta u+\frac{ 2a }{|x|^2}x\cdot\nabla u-\frac{\big(N-2-2a\big)^2}{4}\frac{u}{|x|^{2}}\Big),\quad\forall\, u\in C^2_c(\R^N),
\end{eqnarray*}
where $a<\frac{N-2}{2}$.  If we take  $\mu_1,\mu_2$ to replace  $2a$ and $-\frac14 \big(N-2-2a\big)^2$ respectively  in $\cL_a$,  we propose  a  degenerate operator
 \begin{equation}\label{CNK-0}
 {\cL}_{\mu_1,\mu_2}= -\Delta   +\frac{\mu_1 }{|x|^2}x\cdot\nabla +\frac{\mu_2 }{|x|^2} 
 \end{equation}
with  two free parameters $\mu_1,\mu_2$. Note that  the operator ${\cL}_{\mu_1,\mu_2}$ degenerates at the origin both for the gradient term and the critical Hardy term,  and here we call it {\it the CKN operator}. Our aim is to consider the qualitative properties of the solutions of the elliptic  equations  with that operator.

  \smallskip

 When $\mu_2=0$,  ${\cL}_{\mu_1,0}= -\Delta   +\frac{\mu_1 }{|x|^2}x\cdot\nabla$ is a type of degenerate elliptic operator, which together with  its divergence form,  plays an important role in the harmonic analysis, see \cite{M}, it attracts lots of attentions,  the basic regularities  of related equations  \cite{FKS}, qualitative properties for equation with more general degenerate operators in divergence  form  \cite{STV}.\smallskip

When $\mu_1=0$, ${\cL}_{0,\mu_2}$ reduces to the Hardy-Leray operator, here we can write   ${\cL}_{\mu}=-\Delta+\frac{\mu}{|x|^2}$, which is the prototype of the degenerate operators. The equations with Hardy-Leray operators has been studied extensively in the last decades.   The authors in \cite{GuVe} initiated the analysis of isolated singular solutions of  semilinear problems   $u\mapsto \cL_\mu u+g(u)$  under the condition $\mu\geq -\frac{(N-2)^2}{4}$,  where $\frac{(N-2)^2}{4}$ is the best constant of Hardy inequality (\ref{CNK-in1}),   more related Hardy inequalities refer to \cite{BM,GR,VZ}. Normally, the distributional solution of the Hardy problem $\mathcal{L}_\mu u+g=0\ {\rm in}\  \Omega$
would be  proposed as
$$
 \int_\Omega u\mathcal{L}_\mu \xi\, dx+\int_\Omega g\xi\, dx=0,\quad\forall\,\xi\in C^\infty_c(\Omega),
$$
where  $g$ is a nonlinearity of $x$ and $u$. In this distributional sense, \cite{BP,D,FM1} show the existence   for particular nonlinearity of $u$ under some restriction that  $N\ge3$ and $\mu\in[\mu_0,0)$.    Later on,  C\^irstea at el  in \cite{CC}, C\^irstea in \cite{C}  classified the isolated singular classical solution of $\mathcal{L}_\mu u+b(x)h(u)=0$ in $\Omega\setminus \{0\}$, where
both $b$ and $h$ consist of regularly varying and slowly varying parts (see their definitions in \S1.2.2 of \cite{C}).  There a solution is considered as a $C^1(\Omega \setminus \{0\})$-solution in the sense of distributions in $\Omega \setminus \{0\}$, that is,
 $$
 \int_\Omega \nabla u \nabla \varphi dx - \int_\Omega \frac{\lambda}{|x|^2} u \varphi dx + \int_\Omega b(x)h(u) \varphi dx =0, \quad\forall\,\varphi \in C^1_c(\Omega\setminus \{0\}).
$$
Thanks to a notion of weak solutions of $\cL_\mu u=0$ combined with a  dual formulation of the equation introduced in \cite{CQZ}  the equation
 $$
\cL_\mu u+ g(u)=\nu \quad {\rm in }\, \, \Omega \qquad u=0 \quad {\rm on }\, \partial\Omega,
$$
where $\Omega$ is a bounded smooth domain, $g$ is a continuous nondecreasing function and $\nu$ is a Radon measure. When the pole  of the Leray-Hardy potential is addressed on the boundary of the domain  $\Omega$,  \cite{ChVe1,ChVe2}    extend the approach to classify the boundary isolated singular solutions of Poisson problem
$$
\cL_\mu u+ g(u)=\nu_1 \quad {\rm in }\ \, \Omega,  \qquad u=\nu_2 \quad {\rm on }\ \partial\Omega,
$$
where $\nu_1,\,\nu_2$   are bounded Radon measures respectively on $\Omega$ and $\partial\Omega$.  Recently, the singularities of  the Hardy problems have been also studied extensively  \cite{BMM,BMN,MN,KP,KP1,KV} and the references therein.  \smallskip

When $\mu_1,\mu_2\not=0$, the CKN operator  $ {\cL}_{\mu_1,\mu_2}$ defined in (\ref{CNK-0}) degenerates at the origin thanks to both the gradient term and the Hardy potential. It is worth noting that  the gradient term has the singularity $|x|^{-1}$ indeed, which is also critical at the origin, compared with the Laplacian operator.   The CNK operator in the diverging form $\cL_a$ has been discussed in \cite{FKS,HKM}.  So our aim is to characterize the roles of the two critical  terms in the  related fundamental solution  and the related Poisson problems.

To consider the fundamental solution of $\cL_{\mu_1,\mu_2}$, we provide the following
 setting of $\mu_1,\mu_2$  in this article is the following:
$$\mu_1\in\R,\quad\ \ \mu_2\geq-  \frac{(2-N+\mu_{1})^2}{4}.$$
In this setting, direct computation shows that    the homogeneous problem
\begin{equation}\label{CNK-hom}
\cL_{\mu_1,\mu_2} u=0\quad{\rm in}\ \, \R^N\setminus\{0\}
\end{equation}
has two the radially symmetric solutions
\begin{equation}\label{fun 1}
\Phi_{\mu_1,\mu_2}(x)=\left\{ \arraycolsep=1pt
\begin{array}{lll}
|x|^{\tau_{-}(\mu_1,\mu_2)}\ \ \ \ &{\rm if}\ \, \mu_2>-  \frac{(2-N+\mu_{1})^2}{4},
\\[2mm]
-|x|^{\tau_-(\mu_1,\mu_2)}{\ln|x|}\ \  &{\rm if}\ \, \mu_2=-  \frac{(2-N+\mu_{1})^2}{4}
\end{array}
\right.
\end{equation}
and
\begin{equation}\label{fun 2}
\Gamma_{\mu_1,\mu_2}(x)=|x|^{\tau_{+}(\mu_1,\mu_2)},
\end{equation}
where
$$
 \tau_{-}(\mu_1,\mu_2)=\frac{(2-N+\mu_{1})-\sqrt{(2-N+\mu_{1})^2+4\mu_{2}}}{2}
$$
 and
$$
\tau_{+}(\mu_1,\mu_2)=\frac{(2-N+\mu_{1})+\sqrt{(2-N+\mu_{1})^2+4\mu_{2}}}{2}.
$$
When $\mu_2=-  \frac{(2-N+\mu_{1})^2}{4}$, we observe that
$$ \tau_{-}(\mu_1,\mu_2)= \tau_{+}(\mu_1,\mu_2)=\frac{2-N+\mu_{1}}2:=\tau_0(\mu_1)$$
and
$$\tau_0(\mu_1)<0\ \iff \  \mu_{1}<N-2.$$
In the following, we use the notation of $\tau_{\pm}$ replaced by $\tau_{\pm}(\mu_1,\mu_2)$, $\tau_0$ done by $\tau_0(\mu_1)$ for simplicity if there is no confusion.

\begin{theorem}\label{teo 1}
Assume that $N\geq 2$, $\mu_1\in\R$ and
$$ \mu_2\geq -\frac{(2-N+\mu_{1})^2}{4}.  $$

Let $d\gamma_{\mu_1,\mu_2}:=|x|^{\tau_+-\mu_1 }dx$ and
$$
{\cL}_{\mu_1,\mu_2}^*= -\Delta +(-2\tau_{+}+\mu_{1})\frac{x}{|x|^2}\cdot\nabla,
$$
then we have
\begin{equation}\label{identity 1}
 \int_{{R}^N } \Phi_{\mu_1,\mu_2} \, \cL^*_{\mu_1,\mu_2}(\xi)d\gamma_{\mu_1,\mu_2}=c_{\mu_1,\mu_2}\xi(0),\ \ \ \forall\,\xi\in C^{2}_c(\R^N),
\end{equation}
 where
$$c_{\mu_1,\mu_2}
=\left\{ \arraycolsep=1pt
\begin{array}{lll}
\sqrt{(2-N+\mu_{1})^2+4\mu_{2}}\, |\mathbb{S}^{N-1}|  \ \ \ \ &{\rm if}\ \,  \mu_2> -\frac{(2-N+\mu_{1})^2}{4},
\\[2mm]
 |\mathbb{S}^{N-1}|  \ \  &{\rm if}\ \, \mu_2=-  \frac{(2-N+\mu_{1})^2}{4}.
\end{array}
\right.
 $$

\end{theorem}

\begin{remark}\label{rem 1.1}
 $(i)$  The identity (\ref{identity 1}) means that
$$\cL_{\mu_1,\mu_2}\Phi_{\mu_1,\mu_2}=c_{\mu_1,\mu_2}\delta_0$$
in the weighted distributional sense (\ref{identity 1}),
where $\delta_0$ is Dirac mass concentrates at the origin.\smallskip

$(ii)$ When $\mu_1=0$, the fundamental solution  expressed by Dirac mass in the weighted distributional is derived in \cite{CQZ} and  singular solution of Lane-Emden equations with Hardy operators is classified in \cite{CZ} in this framework.   \smallskip

$(iii)$ When $\mu_2=0$, $\mu_1\leq N-2$, $\Gamma_{\mu_1,\mu_2}=1$ and $\Phi_{\mu_1,\mu_2}(x)=-|x|^{\tau_-}\ln|x|$, which is the fundamental solution of $\cL_{\mu_1,0}$
in the weighted distributional identity (\ref{identity 1}), which coincides the normal distributional identity and $\cL^*_{\mu_1,0}=\cL_{\mu_1,0}$.

\end{remark}

Our second concern is to consider  the  nonexistence  of positive solution for Lane-Emden equation with  CKN operators
\begin{equation}\label{eq 1.1-in}
\left\{\arraycolsep=1pt
\begin{array}{lll}
 \cL_{\mu_1,\mu_2} u\geq Q    u^p \quad
   &{\rm in}\ \ \Omega\setminus\{0\}, \\[2mm]
 \phantom{ \cL_{\mu_1,\mu_2}   }
u\geq0 \quad &{\rm in}\ \ \, \partial\Omega,
 \end{array}
 \right.
\end{equation}
where  $p>0$,   $\Omega$ is a bounded domain containing the origin
  and the potential $Q\in C^\beta_{loc}(\R^N\setminus\{0\})$ with $\beta\in(0,1)$ is a positive function such that for some $\theta>-2$,
\begin{equation}\label{pt r1}
  \liminf_{|x|\to0^+}Q(x)|x|^{-\theta}>0.
\end{equation}

Before stating our main results,   let us involve two  critical exponents
\begin{equation}\label{eq 1.1-cr1}
p^\#_{\mu_1,\mu_2,\theta}=
  1+\frac{2+\theta}{-\tau_+}  \quad {\rm for }\;\  -\frac{(N-2-\mu_1)^2}{4}\leq \mu_2<0.\end{equation}
Here  $p^\#_{\mu_1,\mu_2,\theta}$
is a particular exponent appearing only for $\tau_+<0$, which is the essence for the nonexistence of positive solutions to (\ref{eq 1.1-in}). Note that  any positive solution of
(\ref{eq 1.1-in})  blows up at least like $\Gamma_\mu$ at the origin and this singularity would be improved by the interact of  the nonlinearity $Qu^p$,  which may lead to  an unadmissible singularity in some weighted $L^1$ space.  Inspired by this observation,  \cite{BDT,D} showed the nonexistence of positive solution of  (\ref{eq 1.1-in})   for $p\geq p^*_{0,\mu_2,0}$ when $\mu_1=0$, $\mu_2\in[-\frac{(N-2)^2}{4},0)$  and $Q\equiv 1$.

 Our interest in this paper is  to obtain the nonexistence of positive  solutions of (\ref{eq 1.1-in}) with the  parameters of $\mu_1,\mu_2$ in some suitable range. Here a function  $(u)$ is a positive   solution of   (\ref{eq 1.1-in}), if $u$  satisfies  the inequalities
  $ \mathcal{L}_{\mu_1,\mu_2}  u(x)\geq  Q(x)u^p(x)  \quad  {\rm for\ any \ } x\in  \Omega\setminus\{0\}.$


\begin{theorem}\label{teo 1-non}
Let  $N\geq3$,
$$\mu_1<N-2,\quad \quad   -\frac{(N-2-\mu_1)^2}{4}\leq \mu_2<0,$$
 potential $Q$ verify (\ref{pt r1})  with $\theta>-2$ and $p^\#_{\mu_1,\mu_2,\theta}$ be defined in (\ref{eq 1.1-cr1}).\smallskip

Then for  $p\geq  p^\#_{\mu_1,\mu_2,\theta}$, problem (\ref{eq 1.1-in})    has no positive solution.
\end{theorem}

This type of nonexistence is based on the classification of Poisson problem with the CKN operator
$$\cL_{\mu_1,\mu_2} u=g\quad{\rm in}\ \,  \Omega\setminus\{0\},\qquad u=0\quad{\rm on}\ \,  \partial\Omega,$$
 where $\Omega$ is a bounded domain containing the origin.  We provide sharp conditions of $g$ for the existence and nonexistence of a positive solution to the Poisson problem. The precise results see Section 3 below.
\smallskip



The rest of this paper is organized as follows. In Section 2, we build a weighted distributional identity for $\Phi_{\mu_1,\mu_2}$.    Section 3 is devoted to classify the isolated singular solution of  the related Poisson problem and some important estimates.  Finally,   we deal with the nonexistence of super solutions for semilinear problem (\ref{eq 1.1-in}) in Section 4.

\setcounter{equation}{0}
\section{ Fundamental solutions  }

Here we first remark that the derivation of $\tau_{\pm} =\frac{(2-N+\mu_{1})\pm \sqrt{(2-N+\mu_{1})^2+4\mu_{2}}}{2}$ is based on the calculation
 $$\cL_{\mu_1,\mu_2} |x|^\tau =\Big(-\tau^{2}+(2-N+\mu_{1})\tau+\mu_{2}\Big){|x|^{\tau-2}}$$
and $\tau_{\pm}$ are zero points of
$
 \tau(N-2-\mu_1+\tau)-\mu_{2}=0
$
and $\tau_++\tau_-=2-N+\mu_1$.    Furthermore, we observe that
\begin{eqnarray*}
\cL_{\mu_1,\mu_2} \left(|x|^\tau(-\ln|x|)\right) &=&\Big(-\tau^{2}+(2-N+\mu_{1})\tau+\mu_{2}\Big) {|x|^{\tau-2}(-\ln|x|)}
\\[2mm]&&-\Big(N-2+\mu_1+2\tau\Big)|x|^{\tau-2},
\end{eqnarray*}
which implies that for $ \mu_2=-  \frac{(2-N+\mu_{1})^2}{4} $ and $\tau_0=\frac{2-N+\mu_1}{2}$,
$$\cL_{\mu_1,\mu_2} \left(|x|^\tau(-\ln|x|)\right)=0 \quad{\rm in}\ \, \R^N\setminus\{0\}.$$

From the definition $\cL^*_{\mu_1,\mu_2}$, direct computation shows that
\begin{equation}\label{test bdd}
|\cL^*_{\mu_1,\mu_2}\xi(x)|\leq c\Big(\|\xi\|_{C^2} +\frac1{|x|}\|\xi\|_{C^1}\Big),\quad\forall\, \xi\in C^2_c(\R^N)
\end{equation}
for some $c>0$ independent of $\xi$.

\medskip

\noindent {\bf Proof of Theorem \ref{teo 1}.}
{\it Case 1: $ \mu_2> -\frac{(2-N+\mu_{1})^2}{4}$.}
  Let $u_0=\Phi_{\mu_1,\mu_2}$, then we have that $ \xi \in C^{2}(\R^N)$
\begin{eqnarray*}
0&=&\int_{{\R}^N\setminus B_\varepsilon}  \big({|x|^{\tau_+-\mu_1 }}{\xi})\cL_{\mu_1\mu_2}u_0 dx
 \\[2mm] &=& \int_{{\R}^N\setminus B_\varepsilon}(-\Delta)u_0\, (|x|^{\tau_+-\mu_1 }{\xi)}dx+\int_{{\R}^N\setminus B_\varepsilon}\frac{\mu_1}{|x|^2}{x}\cdot\nabla u_0  (|x|^{\tau_{\tau+}-\mu_1 }{\xi})dx
\\[2mm] &&  +\int_{{\R}^N\setminus B_\varepsilon}\frac{\mu_2}{|x|^2} |x|^{\tau_+-\mu_1 }\xi u_0dx,
 \end{eqnarray*}
where $B_r(z)$ is the ball with radius $r>0$ centered at $z$, particularly, $B_r=B_r(0)$.

 We note that
\begin{eqnarray*}
\mu_1\int_{\R^N\setminus B_\varepsilon}\frac{1}{|x|^2}{x}\cdot\nabla u_0  \,  |x|^{\tau_{+} -\mu_1} {\xi}dx&=&-\mu_1 [N+(\tau_{+}-\mu_{1}-2)]\int_{{\R}^N\setminus B_\varepsilon}|x|^{-N}\xi dx\\[2mm] &&-\mu_1\int_{{\R}^N\setminus B_\varepsilon}|x|^{-N}x\cdot\nabla\xi dx
 -\mu_1\int_{\partial B_\varepsilon} |x|^{1-N }\xi d\omega  \end{eqnarray*}
 and
 $$\int_{{\R}^N\setminus B_\varepsilon}\frac{\mu_2}{|x|^2} |x|^{\tau_+-\mu_1 }\xi u_0dx=\mu_2\int_{{\R}^N\setminus B_\varepsilon} |x|^{-N } \xi dx.  $$
 Direct computation shows
 \begin{eqnarray*}
 \int_{{\R}^N\setminus B_\varepsilon}(-\Delta)u_0\, (|x|^{\tau_+-\mu_1} {\xi})dx
 &=&\int_{{\R}^N\setminus B_\varepsilon}\nabla{u_0}\nabla(|x|^{\tau_+ -\mu_1}\xi)dx+\int_{\partial B_\varepsilon} \frac{x\cdot \nabla{u_0}}{|x|}|x|^{\tau_+-\mu_1 }\xi d\omega
 \\[2mm]&=&\int_{{\R}^N\setminus B_\varepsilon}{u_0}(-\Delta)(|x|^{\tau_+-\mu_1 }\xi)dx+\int_{\partial B_\varepsilon}\frac{x\cdot \nabla{u_0}}{|x|}|x|^{\tau_+-\mu_1 }  \xi  d\omega
 \\[2mm]&&-\int_{\partial B_\varepsilon}\frac{x\cdot\nabla(|x|^{\tau_+ -\mu_1}\xi)}{|x|} u_0d\omega,
  \end{eqnarray*}
 where $\nu=-\frac{x}{|x|}$ be the unit outside normal vector of ${\R}^N\setminus\ B_\varepsilon$ and $d\omega$ be the Hausdorff measure of $\mathbb{S}^{N-1}$,
  \begin{eqnarray*}
 \int_{{\R}^N\setminus B_\varepsilon}{u_0}(-\Delta)(|x|^{\tau_+-\mu_1 }\xi)dx&=& \int_{{\R}^N\setminus B_\varepsilon} u_0  (-\Delta) \xi d\gamma_{\mu_1,\mu_2}-2 (\tau_+-\mu_1) \int_{{\R}^N\setminus B_\varepsilon} u_0   \frac{x}{|x|^2}\cdot\nabla\xi d\gamma_{\mu_1,\mu_2}
 \\[2mm]&&- (\tau_+-\mu_1) (N-2+\tau_+-\mu_1) \int_{{\R}^N\setminus B_\varepsilon}     |x|^{-N } \xi dx,
  \end{eqnarray*}
   \begin{eqnarray*}
   \int_{\partial B_\varepsilon}\frac{x\cdot \nabla{u_0}}{|x|}|x|^{\tau_+-\mu_1 }  \xi  d\omega&=&\tau_-\int_{\partial B_\varepsilon} |x|^{1-N }  \xi  d\omega=\tau_-\int_{\partial B_\varepsilon} |x|^{1-N}  \big(\xi(0)+\nabla \xi(0)\cdot x\big)  d\omega
   \\&\to&\tau_- |\bS^{N-1}|\xi(0)\quad {\rm as}\ \ \varepsilon\to0^+
   \end{eqnarray*}
 and
  \begin{eqnarray*}
  \int_{\partial B_\varepsilon}\frac{x\cdot\nabla(|x|^{\tau_+ -\mu_1}\xi)}{|x|} u_0d\omega&=&(\tau_+-\mu_1) \int_{\partial B_\varepsilon} |x|^{1-N }  \xi  d\omega+\int_{\partial B_\varepsilon} |x|^{1-N } x\cdot \nabla \xi d\omega
  \\&=&(\tau_+-\mu_1)\int_{\partial B_\varepsilon} |x|^{1-N}  \big(\xi(0)+\nabla \xi(0)\cdot x\big)  d\omega+\int_{\partial B_\varepsilon} |x|^{1-N} x\cdot \nabla \xi d\omega
   \\&\to&(\tau_+-\mu_1)  |\bS^{N-1}|\xi(0)\quad {\rm as}\ \ \varepsilon\to0^+,
   \end{eqnarray*}
   where $|\nabla \zeta|$ is bounded.

Now we conclude that
    \begin{eqnarray*}
 0&=& \int_{{\R}^N\setminus B_\varepsilon}(-\Delta)u_0\, (|x|^{ \tau_+-\mu_1 }{\xi)}dx+\int_{{\R}^N\setminus B_\varepsilon}\frac{\mu_1}{|x|^2}{x}\cdot\nabla u_0  (|x|^{ \tau_+-\mu_1}{\xi})dx
 \\&&\ +\int_{{\R}^N\setminus B_\varepsilon}\frac{\mu_2}{|x|^2} |x|^{ \tau_+-\mu_1}\xi u_0dx
  \\&=&\int_{{\R}^N\setminus B_\varepsilon} u_0  (-\Delta) \xi d\gamma_{\mu_1,\mu_2}+(-2\tau_{+}+\mu_{1}) \int_{{\R}^N\setminus B_\varepsilon} u_0   \frac{x}{|x|^2}\cdot\nabla\xi d\gamma_{\mu_1,\mu_2}
 \\&&- \Big(\mu_1 \big(N+(\tau_{+}-\mu_{1}-2)\big)+ (\tau_+-\mu_1) (N-2+\tau_+-\mu_1) -\mu_2\Big)\int_{\R^N\setminus B_\varepsilon} |x|^{-N} \xi dx
  \\&&-(\tau_+-\tau_-)\int_{\partial B_\varepsilon} |x|^{1-N}  \big(\xi(0)+ O(|x|)\big)  d\omega-\int_{\partial B_\varepsilon} |x|^{1-N} x\cdot \nabla \xi d\omega
  \\&=&\int_{{\R}^N\setminus B_\varepsilon } u_0 \cL^*_{\mu_1,\mu_2} \xi d\gamma_{\mu_1,\mu_2}-(\tau_+-\tau_-)\int_{\partial B_\varepsilon} |x|^{1-N}  \big(\xi(0)+ O(|x|)\big)  d\omega
  \\&&- \Big(\tau_+  (N-2+\tau_+-\mu_1) -\mu_2\Big)\int_{\R^N\setminus B_\varepsilon} |x|^{-N} \xi dx
  -\int_{\partial B_\varepsilon} |x|^{1-N} x\cdot \nabla \xi d\omega
   \\&\to&\int_{{\R}^N } u_0 \cL^*_{\mu_1,\mu_2} \xi d\gamma_{\mu_1,\mu_2} -(\tau_+-\tau_-)  |\bS^{N-1}|\xi(0)\quad {\rm as}\ \ \varepsilon\to0^+,
   \end{eqnarray*}
where $\tau_+  (N-2+\tau_+-\mu_1) -\mu_2=0$,
$$\tau_+-\tau_-=\sqrt{(2-N+\mu_{1})^2+4\mu_{2}}$$
and  thanks to (\ref{test bdd})
 $$\int_{{\R}^N\setminus B_\varepsilon } u_0 \cL^*_{\mu_1,\mu_2} \xi d\gamma_{\mu_1,\mu_2}\to \int_{{\R}^N } u_0 \cL^*_{\mu_1,\mu_2} \xi d\gamma_{\mu_1,\mu_2}\quad {\rm as}\ \,  \epsilon\to0^+.$$
Here we note that
$$\sqrt{(2-N+\mu_{1})^2+4\mu_{2}}>0\   \iff\    \mu_2>-\frac{(2-N+\mu_{1})^2}{4}.$$
As a consequence, we obtain that
\begin{equation}\label{identity 2.1}\int_{{\R}^N } u_0 \cL^*_{\mu_1,\mu_2} \xi d\gamma_{\mu_1,\mu_2} =c_{\mu_1,\mu_2}\xi(0),
\end{equation}
 where $$c_{\mu_1,\mu_2}= \Big(\sqrt{(2-N+\mu_{1})^2+4\mu_{2}}\Big){|\mathbb{S}^{N-1}|}.$$

{\it Case 2: $ \mu_2= -\frac{(2-N+\mu_{1})^2}{4}$.}
When
$$\mu_1\in\R \ \ \ \ {\rm and} \ \,  \mu_2= -\frac{(N-2)^2-2(N-2)\mu_1}{4},$$
the   solution $u_0=\Phi_{\mu_1,\mu_2}$ has different form
and then
 \begin{eqnarray*}
   \int_{\partial B_\varepsilon}\frac{x\cdot \nabla{u_0}}{|x|}|x|^{\tau_+-\mu_1 }  \xi  d\omega&=&\tau_-\int_{\partial B_\varepsilon} |x|^{1-N }(-\ln |x|)  \xi  d\omega- \int_{\partial B_\varepsilon} |x|^{1-N }  \xi  d\omega
   \end{eqnarray*}
   and
   \begin{eqnarray*}
  \int_{\partial B_\varepsilon}\frac{x\cdot\nabla(|x|^{\tau_+ -\mu_1}\xi)}{|x|} u_0d\omega&=&(\tau_+-\mu_1) \int_{\partial B_\varepsilon} |x|^{1-N }(-\ln|x|)  \xi  d\omega
  \\&&+\int_{\partial B_\varepsilon} |x|^{1-N }(-\ln|x|) x\cdot \nabla \xi d\omega,
   \end{eqnarray*}
   then
\begin{eqnarray*}
 0&=& \int_{{\R}^N\setminus B_\varepsilon}(-\Delta)u_0\, (|x|^{ \tau_+-\mu_1}{\xi)}dx+\int_{{\R}^N\setminus B_\varepsilon}\frac{\mu_1}{|x|^2}{x}\cdot\nabla u_0  (|x|^{ \tau_+-\mu_1}{\xi})dx
 \\&&\ +\int_{{\R}^N\setminus B_\varepsilon}\frac{\mu_2}{|x|^2} |x|^{ \tau_+-\mu_1}\xi u_0dx
  \\&=&\int_{{\R}^N\setminus B_\varepsilon } u_0 \cL^*_{\mu_1,\mu_2} \xi d\gamma_{\mu_1,\mu_2}-(\tau_+-\tau_-)\int_{\partial B_\varepsilon}|x|^{1-N}(-ln|x|)\xi d\omega- \int_{\partial B_\varepsilon} |x|^{1-N}  \big(\xi(0)+ O(|x|)\big)  d\omega\\&&-\int_{\partial B_\varepsilon} |x|^{1-N}(-\ln|x|) x\cdot \nabla \xi d\omega
   \\&\to&\int_{{\R}^N } u_0 \cL^*_{\mu_1,\mu_2} \xi d\gamma_{\mu_1,\mu_2} -|\bS^{N-1}|\xi(0)\quad {\rm as}\ \ \varepsilon\to0^+,
   \end{eqnarray*}
   where
    $$c_{\mu_1,\mu_2}= {|\mathbb{S}^{N-1}|}
 \quad \ {\rm if} \ \  \mu_2= \frac{\mu_1^2}4-\frac{(2-N+\mu_{1})^2}{4}.$$
We complete the proof. \hfill$\Box$

\setcounter{equation}{0}
\section{  Basic tools }

An important tool is the comparison principle.

\begin{lemma}\label{cr hp}
Assume that $\mu_1\in\R,$ $\mu_2\geq -\frac{(N-2-\mu_1)^2}{4}$, $O$ is a bounded $C^2$ domain containing the origin and $u_i$ with $i=1,2$ are classical solutions of
\begin{equation}\label{eq0 2.1}
 \arraycolsep=1pt\left\{
\begin{array}{lll}
 \displaystyle \mathcal{L}_{\mu_1,\mu_2} u_i = f_i\qquad
   &{\rm in}\quad  O\setminus \{0\},\\[1.5mm]
 \phantom{ L_\mu     }
 \displaystyle  u_i= 0\qquad  &{\rm   on}\ \ \partial O
 \end{array}\right.
\end{equation}
and
$$\limsup_{|x|\to0^+}u_1(x)\Phi_{\mu_1,\mu_2}^{-1}(x)\leq \liminf_{|x|\to0^+}u_2(x)\Phi_{\mu_1,\mu_2}^{-1}(x).$$
  If $f_1\le f_2$ in $O\setminus \{0\}$, then
$$u_1\le u_2\quad{\rm in}\quad O\setminus \{0\}.$$

\end{lemma}
{\bf Proof.}   {\it Case 1: $\mu_2> -\frac{(N-2-\mu_1)^2}{4}$. }   Let $u=u_1-u_2$ satisfy that
$$\mathcal{L}_{\mu_1,\mu_2} u \le 0\quad {\rm in}\ \, O\setminus\{0\}\qquad {\rm and}\quad \limsup_{|x|\to0^+}u(x)\Phi_{\mu_1,\mu_2}^{-1}(x)\leq0,$$
 then for any $\epsilon>0$, there exists $r_\epsilon>0$ converging to zero as $\epsilon\to0$ such that
 $$u\le \epsilon \Phi_{\mu_1,\mu_2}\quad{\rm in}\quad \overline{B_{r_\epsilon}}\setminus\{0\}.$$
We see that
$$u=0<\epsilon \Phi_{\mu_1,\mu_2} \quad{\rm on}\ \ \partial O,$$
then by standard comparison principle for $\mu\geq\mu_0$, see CKN inequality, we have that
 $u\leq \epsilon \Phi_{\mu_1,\mu_2}  $ in $O\setminus\{0\}. $
By the arbitrary of $\epsilon$, we have that $u\leq 0$  in $O\setminus\{0\}$.\smallskip

 {\it Case 2: $\mu_2=-\frac{(N-2-\mu_1)^2}{4}$. }  Replace   $\Phi_{\mu_1,\mu_2}$ in case 1
 by
 $$w_{t_0}:=\Phi_{\mu_1,\mu_2}+t_0\Gamma_{\mu_1,\mu_2},$$
 where $t_0\geq0$ be  such that
 $$w_{t_0}\geq 0\quad {on}\ \partial O$$
 thanks to the boundedness of $O$.

 Repeat the above argument of the case 1 replacing
$\Phi_{\mu_1,\mu_2}$ by $w_{t_0}$.
  \hfill$\Box$\medskip

\subsection{ Poisson problems  }

Our second purpose in this article is to classify the isolated singular solutions of Poisson problem
with CKN operator
\begin{equation}\label{eq p-CKN}
 \left\{ \arraycolsep=1pt
\begin{array}{lll}
 {\cL}_{\mu_1,\mu_2}u=f  \quad{\rm in}\ \, \, \Omega\setminus\{0\},
\\[2mm]
\qquad\ \ u=0 \quad{\rm on}\ \, \partial\Omega,
\end{array}
\right.
\end{equation}
where $f:\Omega\mapsto \R$ is a measurable function and $\Omega$ is a bounded smooth
domain contains the origin. For Poisson problem (\ref{eq p-CKN}), we have the following results.

\begin{theorem}\label{teo 3}
Let  $\mu_1\in\R $ and
$$\mu_2\geq -\frac{(N-2-\mu_1)^2}4 $$
and  $f$ be a function in $C^\gamma_{loc}(\overline{\Omega}\setminus \{0\})$ for some $\gamma\in(0,1)$.

$(i)$ Assume that
\begin{equation}\label{f1}
f \in L^1(\Omega, d\gamma_{\mu_1,\mu_2} )\qquad {\rm i.\,e.}\;\;  \int_{\Omega} |f|\,   d\gamma_{\mu_1,\mu_2}  <+\infty,
\end{equation}
then  for any $k\in\R$  problem (\ref{eq p-CKN}) has a classical solution  $u_k$ with the asymptotic behavior
\begin{equation}\label{singular}
\lim_{x\to0}u_k(x)|x|^{-\tau_-}=k.
\end{equation}

\smallskip

$(ii)$ Assume that $f\ge0$ and
\begin{equation}\label{f2}
 \lim_{r\to0^+} \int_{\Omega\setminus B_r(0)} f\, d \gamma_{\mu_1,\mu_2}  =+\infty,
\end{equation}
then problem (\ref{eq p-CKN}) has no nonnegative solutions.
\end{theorem}

\begin{theorem}\label{teo 4}
  Assume that
$$\mu_1\in\R,\quad \mu_2< -\frac{(N-2-\mu_1)^2}4 $$
and  $f$ be a nonnegative nonzero function,
then problem (\ref{eq p-CKN}) has no nonnegative solutions.
\end{theorem}

Our method is to transform the CKN Poisson problem into a Hardy Poisson problem and the existence and nonexistence could be derived by the results in \cite{CQZ}.

We recall the Poisson problem with Hardy-Leray operator
\begin{equation}\label{eq 1.1f}
 \arraycolsep=1pt\left\{
\begin{array}{lll}
 \displaystyle   \mathcal{L}_\mu u= f\qquad
   {\rm in}\quad  {\Omega}\setminus \{0\},\\[1.5mm]
 \phantom{   L_\mu   }
 \displaystyle  u= 0\qquad  {\rm   on}\quad \partial{\Omega},
 \end{array}\right.
\end{equation}
where $f:\bar\Omega\setminus\{0\}\mapsto \R$ is a H\"{o}lder continuous locally in $\bar\Omega\setminus\{0\}$ and
$$\mathcal{L}_\mu=-\Delta+\frac{\mu}{|x|^2}$$
is the Hardy-Leray operator with $\mu\geq -\frac{(N-2)^2}{4}$. \cite{CQZ} classifies the    isolated singularities  of solutions to (\ref{eq 1.1f})  by building the connection with weak solutions of
\begin{equation}\label{eq 1.2}
 \arraycolsep=1pt\left\{
\begin{array}{lll}
 \displaystyle   \mathcal{L}_\mu u= f+c_\mu k\delta_0\qquad
   &{\rm in}\quad  {\Omega},\\[1.5mm]
 \phantom{   L_\mu   }
 \displaystyle  u= 0\qquad  &{\rm   on}\quad \partial{\Omega}
 \end{array}\right.
\end{equation}
in the $|x|^{\tau_+(\mu)}dx $-distributional sense, that is, $u\in L^1(\Omega, |x|^{\tau_+(\mu)}dx )$ and satisfying
 \begin{equation}\label{1.2f}
 \int_{\Omega}u  \mathcal{L}_\mu^*(\xi)\, |x|^{\tau_+(\mu)}dx   = \int_{\Omega} f  \xi\, |x|^{\tau_+(\mu)}dx  +c_\mu k\xi(0),\quad\forall\, \xi\in   C^{1.1}_0(\Omega),
\end{equation}
where $k\in\R$ and $\mathcal{L}_\mu^* u=-\Delta u -2\tau_+(\mu) \frac{x}{|x|^2}\cdot \nabla u$ and
$$\tau_\pm(\mu)=\frac{2-N}{2}\pm\sqrt{\frac{(N-2)^2}{4}+\mu}.$$
It is worth noting that $\tau_\pm(\mu)=\tau_\pm(0,\mu)$ with $\mu_1=0$
and $\mu_2=\mu$.

The classification of isolated singular solutions of (\ref{eq 1.1f}) states as following:
\begin{theorem}\label{teo h3}\cite[Theorem 1.3]{CQZ}
Let    $f$ be a function in $C^\gamma_{loc}(\overline{\Omega}\setminus \{0\})$ for some $\gamma\in(0,1)$.

$(i)$ Assume that
\begin{equation}\label{f1}
f \in L^1(\Omega, |x|^{\tau_+(\mu)}dx )\qquad {\rm i.\,e.}\;\;  \int_{\Omega} |f(x)|\,    |x|^{\tau_+(\mu)}dx  <+\infty,
\end{equation}
then  for any $k\in\R$, problem (\ref{eq 1.2}) admits a unique weak solution  $u_k$, which is a classical solution of   problem (\ref{eq 1.1f}).
Furthermore, if assume more that
\begin{equation}\label{4.1-3}
 \lim_{|x|\to0^+}f(x)|x|^{2-\tau_-(\mu)}=0,
\end{equation}
then  we have the asymptotic behavior
\begin{equation}\label{singular}
\lim_{|x|\to0^+}u_k(x)|x|^{-\tau_-(\mu)}(x)=k.
\end{equation}

$(ii)$ Assume that $f$ verifies (\ref{f1})  and $u$ is a nonnegative solution of (\ref{eq 1.1f}), then  $u$ is a weak solution of  (\ref{eq 1.2}) for some $k\ge0$.\smallskip

$(iii)$ Assume that $f\ge0$ and
\begin{equation}\label{f2}
 \lim_{r\to0^+} \int_{\Omega\setminus B_r(0)} f(x)\,  |x|^{\tau_+(\mu)}dx  =+\infty,
\end{equation}
then problem (\ref{eq 1.1f}) has no nonnegative solutions.
\end{theorem}

Now we are in a position to show the isolated singular solution of (\ref{eq p-CKN}). \medskip

\noindent{\bf Proof of Theorem \ref{teo 3}. }
Let $\tau= \frac{\mu_1}{2}$ and $u\in C^2(\Omega\setminus\{0\})$ be a solution of (\ref{eq p-CKN}) and
 $$ u(x)=v(x)|x|^{\tau}\quad{\rm in}\ \, \Omega\setminus\{0\},$$
then
we have that
\begin{eqnarray*}
 \cL_{\mu_1,\mu_2}(|x|^{\tau}v(x))&=&|x|^{\tau} \Big(-\Delta v+(\mu_1-2\tau)\frac{x}{|x|^2}x\cdot \nabla v+\big(\mu_2+\tau \mu_1-\tau(N-2+\tau)\big)\frac{v}{|x|^2}\Big)
 \\&=&|x|^{\tau}  {\cL}_{\tilde \mu } v  ,
 \end{eqnarray*}
 where
 $${\cL}_{\tilde \mu } v=-\Delta v+\tilde \mu\frac{v}{|x|^2}$$
 with
$$\tilde \mu=\mu_2+\frac{ \mu_1^2}4-\frac{\mu_1}{2}(N-2). $$
Hence $v$ is a classical solution of
\begin{equation}\label{eq 4.1}
 \left\{ \arraycolsep=1pt
\begin{array}{lll}
 {\cL}_{\tilde \mu }v=|x|^{-\frac{\mu_1}{2}} f  \quad&{\rm in}\ \, \, \Omega\setminus\{0\},
\\[2mm]
\quad\,  u=0 \quad&{\rm on}\ \, \partial\Omega,
\end{array}
\right.
\end{equation}
where we note that
$\tilde \mu\geq -\frac{(N-2)^2}{4}$ is equivalent
to
$$\mu_2\geq -\frac{(N-2-\mu_1)^2}4. $$
It is remarkable that
\begin{eqnarray*}
\tau_\pm(\tilde \mu)&=&-\frac{N-2}{2}\pm\sqrt{\frac{(N-2)^2}{4}+\mu_2+\frac{ \mu_1^2}4-\frac{\mu_1}{2}(N-2) }
\\&=&\frac{2-N\pm\sqrt{(2-N+\mu_{1})^2+4\mu_{2}}}{2}
\\&=&\tau_\pm(\mu_1,\mu_2)+\frac{\mu_1}2,
\end{eqnarray*}
which implies that
$$f \in L^1(\Omega, d\gamma_{\mu_1,\mu_2} )\ \iff\ |x|^{-\frac{\mu_1}{2}}f \in L^1(\Omega, |x|^{\tau_+(\tilde \mu)}dx )$$
and for $k\in\R$
$$
\lim_{|x|\to0^+}v(x)|x|^{-\tau_-(\tilde\mu)}(x)=k
\ \iff \ \lim_{|x|\to0^+}u(x)|x|^{-\tau_-(\tilde\mu)-\frac{\mu_1}{2}}=\lim_{|x|\to0^+}u(x)|x|^{-\tau_-(\mu_1,\mu_2)-\mu_1 }=k.$$
Therefore, the existence and nonexistence of problem (\ref{eq p-CKN}) follows Theorem \ref{teo h3} part $(i)$ and $(iii)$.
We complete the proof.\hfill$\Box$\medskip

\noindent{\bf Proof of Theorem \ref{teo 4}. } Observe that
$$ \mu_2< -\frac{(N-2-\mu_1)^2}4\ \iff\  \tilde \mu< -\frac{(N-2)^2}{4}$$
and as shown in the proof Theorem \ref{teo 3},
the nonexistence of  problem (\ref{eq p-CKN}) follows by the nonexistence of positive solution
 of (\ref{eq 4.1}) with $\tilde \mu<-\frac{(N-2)^2}{4}$ from \cite[Proposition 5.2]{CQZ}.\hfill$\Box$\medskip

\subsection{Some estimates}

In order to improve the blowing up rate at the origin or decay at infinity, we need the following estimates.   Let
\begin{equation}
  \label{eq:fractional-power}
 {\bf c}(\tau) =-\tau   (N-2 -\mu_1+\tau) +\mu_2  >0,
\end{equation}
the map: $\tau\mapsto c(\tau)$ is concave and
$\tau_+\geq \tau_-$ are the two zero points of  $ {\bf c}(\tau)=0.$
There holds
\begin{equation}
  \label{eq:fractional-power-1}
  \cL_{\mu_1,\mu_2}  |\cdot|^\tau =  {\bf c}(\tau)|\cdot|^{\tau-2 } \qquad  {\rm in}\ \   \R^N\setminus\{0\}.
\end{equation}

   \begin{lemma}\label{lm 3.1-it} Assume that  $\mu_1\in\R,$ $\mu_2\geq -\frac{(N-2-\mu_1)^2}{4}$,
 the nonnegative function $g\in C^\beta_{loc}(B_r \setminus \{0\})$ for some $\beta\in(0,1)$ and there exists $\tau\in(\tau_- ,\tau_+ )$, $c>0$ and $r_0>0$ such that
 $$ g(x)\geq c|x|^{\tau-2}\quad{\rm in}\ \, B_r \setminus \{0\}.$$

Let  $u_g$ be a positive solution of problem
 \begin{equation}\label{eq 2.1-hom}
\cL_{\mu_1,\mu_2} u\geq g \quad {\rm in}\ \ B_r \setminus \{0\},\qquad  u\geq0 \quad  {\rm on}\ \  \partial B_r,
 \end{equation}
then there exists $c_1>0$  such that
$$u_g(x)\geq c_1|x|^\tau\quad{\rm in}\ \  B_{\frac r2}.$$

\end{lemma}
{\bf Proof. }  For $\tau\in(\tau_- ,\tau_+ )$, we have that
$$\cL_{\mu_1,\mu_2} |x|^{\tau}={\bf c}(\tau) |x|^{\tau-2}\quad  {\rm in}\ \  \R^N\setminus\{0\}.$$

 For $r>0$, let
$$w(x)=|x|^{\tau}-r^{\tau-\tau_+}|x|^{\tau_+}\quad  {\rm in}\ \  \R^N\setminus\{0\},$$
which verifies that
$$\cL_{\mu_1,\mu_2} w={\bf c}(\tau)|x|^{\tau-2}\quad{\rm in}\ B_r\setminus\{0\},\qquad
w= 0\ \ {\rm on} \ \ \partial B_r,$$
where ${\bf c}(\tau)>0$.
 Then our argument follows by   Lemma \ref{cr hp}.
 \hfill$\Box$

\setcounter{equation}{0}
\section{ Nonexistence in  a punctured domain }

We prove the nonexistence of positive solution of (\ref{eq 1.1-in}) in a punctured domain $\Omega\setminus\{0\}$ by contradiction, i.e.   (\ref{eq 1.1-in}) is assumed to have a positive solution $u_0$ and we will obtain a contradiction from  Theorem \ref{teo 3} and Theorem \ref{teo 4}.
In this section, we  can assume that  $$B_2 \subset \Omega$$
and
$$Q(x)\geq q_0|x|^{\theta}\quad{\rm in}\ \, B_1 \setminus\{0\}.$$

 \smallskip

\noindent{\bf Proof of Theorem \ref{teo 1}. }
By contradiction, we assume that $u_0$ is a positive super solution of (\ref{eq 1.1-in}) in $\Omega\setminus\{0\}$.
 Let $$v_0(x)=|x|^{\tau_+}-1\quad {\rm for}\ \,  B_1\setminus\{0\},$$
  then
 $$ \cL_{\mu_1,\mu_2} v_0(x)=\frac{\mu_2}{|x|^2}<0 $$
for $\mu_2<0$.  Here we recall that $\tau_+<0$ for $ -\frac{(N-2-\mu_1)^2}{4}\leq \mu_2<0$.

Since $B_2\subset \Omega$ and there exists $t_1>0$ such that
$$u_0\geq t_1\quad {\rm for}\ \ |x|=1.$$
Therefore, for any $\epsilon>0$, there holds
$$\liminf_{x\to0}\big( u_0-t_1v_0+\epsilon \Phi_{\mu_1,\mu_2}\big)(x) \Phi_{\mu_1,\mu_2}^{-1}\geq \epsilon,$$
$$\big( u_0-t_1v_0+\epsilon \Phi_{\mu_1,\mu_2}\big)(x)\geq0\quad {\rm for}\ \, |x|=1$$
and
$$\cL_{\mu_1,\mu_2} ( u_0-t_1v_0+\epsilon \Phi_{\mu_1,\mu_2})\geq0\quad{\rm in} \ B_1\setminus\{0\}$$
It follows by Lemma \ref{cr hp}   that
$$ u_0-t_1v_0+\epsilon \Phi_{\mu_1,\mu_2}\geq0\quad{\rm in} \ B_1\setminus\{0\}$$
and by the arbitrary of $\epsilon>0$, we have that
$$u_0 \geq t_1v_0\quad{\rm in} \ B_1\setminus\{0\}.$$

Let
$$q^\#_{\mu_1,\mu_2,\theta}=  \frac{N+\theta}{-\tau_+}-1.  $$
Then
$$q^\#_{\mu_1,\mu_2,\theta}=p^\#_{\mu_1,\mu_2,\theta}\quad {\rm for}\ \  \mu_2=-\frac{(N-2-\mu_1)^2}{4}$$
and
$$q^\#_{\mu_1,\mu_2,\theta }>p^\#_{\mu_1,\mu_2,\theta}\quad {\rm for}\ \  -\frac{(N-2-\mu_1)^2}{4}< \mu_2<0.$$

Set
$g(x)=Q(x) u_0(x)^{p}$,
then  for some $d_0>0$
$$u_0(x)\geq d_0|x|^{\tau_+}\quad{\rm in}\ \ B_{\frac12} \setminus\{0\}.$$

\noindent{\it Part 1: Nonexistence for $-\frac{(N-2-\mu_1)^2}{4}\leq \mu_2<0$ and $p\geq q^\#_{\mu_1,\mu_2,\theta}$. }
Note that
 \begin{eqnarray*}
 \cL_{\mu_1,\mu_2} u_0(x) \geq  g(x)\geq    d_0^p  |x|^{\theta+\tau_+p}  \quad{\rm in}\ \ B_{\frac12} \setminus\{0\},
 \end{eqnarray*}
 where $\theta+\tau_+p+ \tau_+\leq -N$ and
$$\lim_{r\to0^+}\int_{B_{r_0}(0)\setminus B_r(0)} Q(x)|x|^{\tau_+p}\Gamma_{\mu_1,\mu_2} dx=+\infty$$
by the fact that $p\geq q^\#_{\mu_1,\mu_2,\theta }$.
As a consequence, we see that $u_0$ is a solution of
$$
\cL_{\mu_1,\mu_2} u_0 \geq  g\quad      {\rm in}\ \,  \Omega\setminus\{0\},\qquad u_0\geq0\quad {\rm in}\ \ \R^N\setminus\{0\}.
$$
 Then there is a contradiction from Theorem \ref{teo 3} part $(ii)$. \smallskip

\noindent {\it Part  2: Nonexistence for    $-\frac{(N-2-\mu_1)^2}{4}< \mu_2<0$ and $p\in \big(p^\#_{\mu_1,\mu_2,\theta},\, q^\#_{\mu_1,\mu_2,\theta}\big)$. }
Let $\tau_0=\tau_+<0$, then
 $$
\cL_{\mu_1,\mu_2}  u_0(x)  \geq    q_0d_0^p|x|^{p\tau_0+\theta}=q_0d_0^p |x|^{\tau_1-2}\quad {\rm in}\ \, B_{\frac12}\setminus\{0\},
$$
where
$$\tau_1:=p\tau_0+\theta+2.$$
 By Lemma \ref{lm 3.1-it}, we have that
$$u_0(x)\geq d_1|x|^{\tau_1}\quad {\rm in}\ \, B_{\frac12} \setminus\{0\}.$$
Iteratively, we recall that
$$\tau_j:=p\tau_{j-1}+\theta+2,\quad j=1,2,\cdots .$$
Note that
$$\tau_1-\tau_0=(p-1)\tau_0+\theta+2 <0$$
for $p\in(p^\#_{\theta,\mu}, q^\#_{\theta,\mu})$.

If $\tau_{1}p   +\theta +2 \leq \tau_-,$
then $\tau_{1}p  + \tau_++\theta<-N$ and
\begin{eqnarray*}
 \cL_{\mu_1,\mu_2}  u_0(x) \geq  g(x)\geq   q_0 d_1^p  |x|^{\theta+\tau_{1}p}  \quad{\rm in}\ \ B_{\frac12} \setminus\{0\}
 \end{eqnarray*}
 and  a contradiction comes from Theorem \ref{teo 3} part $(ii)$.

If not, we iterate above procedure. If
$$\tau_{j+1}:=\tau_jp+\theta+2\in  (\tau_-,\tau_+),$$
it following by Lemma \ref{lm 3.1-it}  that
$$u_0(x)\geq d_{j+1} |x|^{\tau_{j+1}},$$
where
$$\tau_{j+1}=p\tau_j+2+\theta<\tau_j.$$
If $\tau_{j+1}p+ \tau_++\theta \leq -N,$
then
\begin{eqnarray*}
 \cL_{\mu_1,\mu_2} u_0(x) \geq    q_0  d_{j+1}^p  |x|^{\theta+\tau_{j+1}p}  \quad{\rm in}\ \ B_{\frac12} \setminus\{0\}
 \end{eqnarray*}
 then we have a contradiction and  we are done.

Furthermore, this iteration must stop by finite times since
$$
\tau_j-\tau_{j-1} = p(\tau_{j-1}-\tau_{j-2})=p^{j-1} (\tau_1-\tau_0)\to-\infty\ \ {\rm as} \ \, j\to+\infty.
$$
As a consequence, we obtain the nonexistence for the case $p\in \big(p^\#_{\mu_1,\mu_2,\theta},\, q^\#_{\mu_1,\mu_2,\theta}\big)$.

 \smallskip

\noindent{\it Part  3:  Nonexistence for    $-\frac{(N-2-\mu_1)^2}{4}< \mu_2<0$ and  $p=p^\#_{\mu_1,\mu_2,\theta}>1$. }
Note that
 \begin{eqnarray*}
 \cL_{\mu_1,\mu_2} u_0(x) \geq       q_0 d_0^p  |x|^{\theta+\tau_+p}  \quad{\rm in}\ \ B_{\frac12} \setminus\{0\},
 \end{eqnarray*}
 where in this case
 $$\theta+\tau_+p+2=\tau_+.$$
Here we see that for some $\sigma_0>0$
 $$\frac12 Qu_0^{p-1}(x)\geq \sigma_0|x|^{-2} \quad{\rm in}\ \ B_{\frac12} \setminus\{0\}.$$
If $\mu_2-\sigma_0< -\frac{(N-2-\mu_1)^2}{4}$,
then we have that
$$\cL_{\mu_1,\mu_2-\sigma_0} u_0(x) =f\geq 0, $$
which has no positive solution by Theorem \ref{teo 4}.

If $\mu_2-\sigma_0\geq -\frac{(N-2-\mu_1)^2}{4}$,   we can write problem (\ref{eq 1.1-in}) as following
\begin{equation}\label{eq 1.1 new}
 \cL_{\mu_1,\mu_2-\sigma_0}  u_0\geq \frac12 Q(x)   u_0^p \quad {\rm in}\ \  B_{\frac12} \setminus\{0\} ,
 \end{equation}
 which the critical exponent
 $$p^\#_{\mu_1,\mu_2-\sigma_0,\theta}=  1+\frac{2+\theta}{-\tau_+(\mu_1,\mu_2-\sigma_0)} <   1+\frac{2+\theta}{-\tau_+(\mu_1,\mu_2 )}=p^\#_{\mu_1,\mu_2,\theta}, $$
 since $t\mapsto \tau_+(\mu_1,t )$ is decreasing.
Thus it reduces to the case:    {\it part 2} for (\ref{eq 1.1 new}) for $p=p^\#_{\mu_1,\mu_2,\theta}$ is supercritical and  we obtain the nonexistence.
\hfill$\Box$\medskip

  \bigskip\bigskip

  \noindent{\bf \small Acknowledgements:}    {\footnotesize
This work is supported by the Natural Science Foundation of China, No. 12071189, by Jiangxi Province Science Fund for Distinguished Young Scholars, No. 20212ACB211005. }

\end{document}